\input amstex\documentstyle{amsppt}  
\pagewidth{12.5cm}\pageheight{19cm}\magnification\magstep1
\topmatter
\title Involutions in Weyl groups and nil-Hecke algebras\endtitle
\author George Lusztig and David A. Vogan, Jr\endauthor
\address{Department of Mathematics, M.I.T., Cambridge, MA 02139}\endaddress
\thanks{GL supported by NSF grant DMS-1855773 and by a Simons Fellowship}\endthanks
\endtopmatter   
\document

\define\part{\partial}

\define\m{\mapsto}
\define\do{\dots}

\define\lra{\leftrightarrow}

\define\sub{\subset}

\define\nl{\newline}
\redefine\i{^{-1}}

\define\dom{\text{\rm dom}}

\define\p{\pi}
\define\ph{\phi}

\redefine\aa{\bold a}

\define\II{\bold I}

\define\NN{\bold N}

\define\TT{\bold T}

\define\ZZ{\bold Z}

\define\ct{\Cal T}

\define\bul{\bullet}

\head Introduction \endhead
\subhead 0.1\endsubhead
Let $W$ be a Coxeter group and let $S$ be the set of simple reflections of $W$;
we assume that $S$ is finite. Let $w\m|w|$ be the length function on $W$.
Let $H$ be the Iwahori-Hecke algebra attached
to $W$. Recall that $H$ is the free $\ZZ[u]$-module with basis $\{\TT_w;w\in W\}$
($u$ is an indeterminate) with (associative) multiplication characterized by
$\TT_s\TT_w=\TT_{sw}$ if $s\in S,w\in W,|sw|=|w|+1$,
$\TT_s\TT_w=u^2\TT_{sw}+(u^2-1)\TT_w$ if $s\in S,w\in W,|sw|=|w|-1$.

Let $w\m w^*$ be an automorphism with square $1$ of $W$
preserving $S$ and let $\II_*=\{w\in W;w^*=w\i\}$ be the set of ``twisted
involutions'' in $W$. Let $M$ be the free $\ZZ[u]$-module with basis
$\{\aa_x;x\in \II_*\}$. 
For any $s\in S$ we define a $\ZZ[u]$-linear map $\TT_s:M@>>>M$ by

$\TT_s\aa_x=u\aa_x+(u+1)\aa_{sx}$ if $x\in\II_*,sx=xs^*,|sx|=|x|+1$,

$\TT_s\aa_x=(u^2-u-1)\aa_x+(u^2-u)\aa_{sx}$ if $x\in\II_*,sx=xs^*,|sx|=|x|-1$,

$\TT_s\aa_x=\aa_{sxs^*}$ if $x\in\II_*,sx\ne xs^*,|sx|=|x|+1$,

$\TT_s\aa_x=(u^2-1)\aa_x+u^2\aa_{sxs^*}$ if $x\in\II_*,sx\ne xs^*,|sx|=|x|-1$.

It is known that the maps $\TT_s$ define an $H$-module structure on $M$. (See
\cite{LV12} for the case where $W$ is a Weyl group or an affine Weyl group and
\cite{L12} for the general case; the case where
$W$ is a Weyl group and $u$ is specialized to $1$
was considered earlier in \cite{K00}.)
When $u$ is specialized to $0$, $H$ becomes the 
free $\ZZ$-module $H_0$ with basis $\{T_w;w\in W\}$ 
with (associative) multiplication characterized by

$T_sT_w=T_{sw}$ if $s\in S,w\in W,|sw|=|w|+1$,

$T_sT_w=-T_w$ if $s\in S,w\in W,|sw|=|w|-1$ 
\nl
(a nil-Hecke algebra).
From these formulas we see that there is a well defined monoid structure
$w,w'\m w\bul w'$ on $W$ such that for any $w,w'$ in $W$ we have

$T_wT_{w'}=(-1)^{|w|+|w'|+|w\bul w'|}T_{w\bul w'}$
\nl
(equality in $H_0$). In this monoid we have

(a) $(w\bul w')\i=w'{}\i\bul w\i$

(b) $(w\bul w')^*=w^*\bul w'{}^*$
\nl
for any $w,w'$ in $W$.

\subhead 0.2\endsubhead
When $u$ is specialized to $0$, the $H$-module $M$ becomes the 
$H_0$-module $M_0$ with $\ZZ$-basis $\{a_w;w\in\II_*\}$ in which for $s\in S$ we
have
$$\align& T_sa_x=a_{sx} \text{ if }x\in\II_*,sx=xs^*, |sx|=|x|+1,\\&
T_sa_x=a_{sxs^*} \text{ if }x\in\II_*,sx\ne xs^*, |sx|=|x|+1,\\&
T_sa_x=-a_x \text{ if }x\in\II_*,|sx|=|x|-1.\tag a\endalign$$
By \cite{L12, 4.5} there is a unique function $\ph:\II_*@>>>\NN$ such that
$\ph(1)=0$ and such that for any $s\in S,x\in\II_*$ with $|sx|=|x|-1$ we have
$\ph(x)=\ph(sx)+1$ if $sx=xs^*$, $\ph(sxs^*)=\ph(x)$ if $sx\ne xs^*$;
moreover we have $|x|=\ph(x)\mod2$ for all $x\in\II_*$.
For $x\in\II_*$ we set $||x||=(1/2)(|x|+\ph(x))\in\NN$. (See also \cite{V79, p.92}.)

One of our main observations is that the action of any
$T_w,w\in W$ on a basis element $a_x,x\in\II_*$ of $M_0$
has a simple description in terms of the monoid $(W,\bul)$,
namely

(b) $T_wa_x=(-1)^{|w|+||x||+||w\bul x\bul w^*{}\i||}
a_{w\bul x\bul w^*{}\i}$.
\nl
(Note by 0.1(a),(b), we have  $w\bul x\bul w^*{}\i\in\II_*$ whenever $x\in \II_*$.) See \S1 for a proof.

\subhead 0.3\endsubhead
In \cite{L16} it is shown that there is a unique map $\p:W@>>>\II_*$ such that
for any $w\in W$ we have $T_wa_1=\pm a_{\p(w)}$ in $M_0$. This can be deduced also
from 0.2(b), which gives a closed formula for $\p$, namely
$$\p(w)=w\bul w^*{}\i.\tag a$$
By \cite{L16, 1.8(c)},

(b) the map $\p:W@>>>\II_*$ is surjective.
\nl
More generally, let $J\sub S$ be such that the subgroup $W_J$ generated by $J$ is finite.
Let $w_J$ be the longest element of $W_J$. Let $J^*$ be the image of $J$ under $*$.
Let ${}^JW=\{w\in W;|w|=|w_J|+|w_Jw|\}$,
$W^{J^*}=\{w\in W;|w|=|w_{J^*}|+|ww_{J^*}|\}$, ${}^JW^{J^*}={}^JW\cap W^{J^*}$. The following extension of (b) is verified in \S2.

(c) $\p$ restricts to a surjective map ${}^J\p:{}^JW@>>>\II_*\cap {}^JW^{J^*}$.

\subhead 0.4\endsubhead
In 2.3 it is shown that when $W$ is an irreducible affine Weyl group then for a suitable
$J,*$, the map in 0.3(c) can be interpreted as a (surjective) map from the set of
translations in $W$ to the set of dominant translations in $W$. This map is bijective
if $W$ is of affine type $A_1$ (see 2.4)
but is not injective if $W$ is of affine type $A_2$ (see 2.5).
This map takes any dominant translation to its square. It would be interesting to find
a simple formula for this map
extending the formulas in 2.4, 2.5.

\head 1. Proof of 0.2(b)\endhead
\subhead \endsubhead
We prove 0.2(b) by induction on $|w|$. If $w=1$ we have
$w\bul x\bul w^*{}\i=x$ hence the desired result holds.
Assume now that $w=s\in S$. If $sx=xs^*,|sx|=|x|+1$, we have
$s\bul x\bul s^*=(sx)\bul s^*=x\bul s^*\bul s^*=x\bul s^*=xs^*=sx$ and
$$\align&1+||x||+||sx||=1+(1/2)(|x|+\ph(x))+1/2(|sx|+\ph(sx))\\&=
1+(1/2)(|x|+\ph(x))+1/2(|x|+1+\ph(x)+1)=|x|+\ph(x)+2=0\mod2.\endalign$$
If $sx\ne xs^*,|sx|=|x|+1$, we have $|sxs^*|=|sx|+1$ hence
$s\bul x\bul s^*=(sx)\bul s^*=sxs^*$ and
$$\align&1+||x||+||sxs^*||=1+(1/2)(|x|+\ph(x))+1/2(|sxs^*|+\ph(sxs^*))\\&=
1+(1/2)(|x|+\ph(x))+1/2(|x|+2+\ph(x))=|x|+\ph(x)+2=0\mod2.\endalign$$
If $|sx|=|x|-1$, we have $|xs^*|=|x|-1$ hence
$s\bul x\bul s^*=x\bul s^*=x$ and
$1+||x||+||x||=1\mod2$ hence the desired result holds.

Assume now that $w\ne1$. We can find $s\in S$ such that
$|sw|=|w|-1$. By the induction hypothesis we have
$$T_{sw}a_x=(-1)^{|sw|+||x||+||(sw)\bul x \bul(sw)^*{}\i||}
a_{(sw)\bul x\bul (sw)^*{}\i}.$$
Using the earlier part of the proof we have
$$\align&
T_wa_w=T_sT_{sw}a_x=(-1)^{|sw|+||x||+||(sw)\bul x\bul(sw)^*{}\i||}
T_sa_{(sw)\bul x\bul(sw)^*{}\i}\\&=
(-1)^{|sw|+||x||+||(sw)\bul x\bul(sw)^*{}\i||}
(-1)^{1+||(sw)\bul x\bul(sw)^*{}\i||+||w\bul x\bul w^*{}\i||}\\&
a_{s\bul(sw)\bul x\bul(sw)^*{}\i\bul s^*}\\&=
(-1)^{|w|+||x||+||w\bul x\bul w^*{}\i||}
a_{w\bul x\bul w^*{}\i}.\endalign$$
This completes the proof of 0.2(b).

\head 2. The map ${}^J\p$\endhead
\subhead 2.1\endsubhead
For $w_1,w_2$ in $W$ we say that $w_1$ is an initial segment of $w_2$ if there exist
$s_1,s_2,\do s_n$ in $S$ and $k\in[0,n]$ such that $w_1=s_1s_2\do s_k$, $w_2=s_1s_2\do s_n$, 
$|w_1|=k,|w_2|=n$; we say that $w_1$ is a final segment of $w_2$ if $w_1\i$ is an initial segment of $w_2\i$. We show:

(a) For $w,w'$ in $W$, $w$ is an initial segment of $w\bul w'$ and
$w'$ is a final segment of $w\bul w'$.
\nl
We argue by induction on $|w'|$. If $w'=1$ the result is obvious. Assume now that $w'\ne1$. We
can find $s\in S$ such that $|w'|=|sw'|+1$. We have $w\bul w'=w\bul s\bul(sw')$.
If $|ws|=|w|+1$ then $w\bul w'=(ws)\bul(sw')$ and by the induction hypothesis $ws$ is an initial
segment of $(ws)\bul(sw')=w\bul w'$. Since $w$ is an initial segment ot $ws$ it follows that
 $w$ is an initial segment of $w\bul w'$. 
If $|ws|=|w|-1$ then $w\bul w'=w\bul(sw')$ and by the induction hypothesis $w$ is an initial
segment of $w\bul(sw')=w\bul w'$. This proves the first assertion of (a). The second assertion of
(a) follows from the first using 0.1(a).

\subhead 2.2\endsubhead
We now fix $J\sub S$ as in 0.3. Let $w\in{}^JW$. Then $w_J$ is an initial segment of $w$
and (by 2.1(a)) $w$ is an initial segment of $w\bul w^*{}\i$ hence
$w_J$ is an initial segment of $w\bul w^*{}\i$ so that $w\bul w^*{}\i\in {}^JW$.
Since $w\bul w^*{}\i\in\II_*$ we see that $w_{J^*}$ is a final segment of $w\bul w^*{}\i$
so that $w\bul w^*{}\i\in W^{J^*}$. Thus we have $w\bul w^*{}\i\in \II_*\cap {}^JW^{J^*}$.
We see that the map ${}^J\p:{}^JW@>>>\II_*\cap {}^JW^{J^*}$ in 0.3(c) is well defined.

We now prove that this map is surjective. Let $x\in \II_*\cap {}^JW^{J^*}$.
Let $z$ be the unique element of minimal length in $W_JxW_{J^*}$.
Now $z^*{}\i$ is again an element of minimal length in $W_JxW_{J^*}$ so it must be equal to $z$.
Thus we have $z\in\II_*$. 
By 0.2(b) we have $T_{w_J}a_z=\pm a_{w_J\bul z \bul w_{J^*}}$.
Note that
$$w_J\bul z \bul w_{J^*}\in W_JxW_{J^*}=W_JzW_{J^*}.$$
By 2.1(a), $w_J$ is an initial segment of $w_J\bul z \bul w_{J^*}$ so that
$$w_J\bul z \bul w_{J^*}\in{}^JW.$$
Since $w_J\bul z \bul w_{J^*}\in\II_*$ we have also
$w_J\bul z \bul w_{J^*}\in W^{J^*}$ so that
$w_J\bul z \bul w_{J^*}\in {}^JW^{J^*}$. Thus
$w_J\bul z \bul w_{J^*}$ is the element of maximal length in $W_JxW_{J^*}$ so that it must be 
equal to $x$ and we have $T_{w_J}a_z=\pm a_x$.
By 0.2(b) we have $T_ea_1=\pm a_z$ for some $e\in W$.
We then have $\pm a_x=\pm T_{w_J}T_ea_1=\pm T_wa_1$
where $w=w_J\bul e$ has $w_J$ as initial segment (see 2.1(a)) so that $w\in {}^JW$.
This proves the surjectivity of ${}^J\p$.

\subhead 2.3\endsubhead
In the remainder of this section we assume that $W$ is an irreducible affine Weyl group.
Let $\ct$ be the (normal) subgroup of $W$ consisting of translations.
We fix a proper subset $J\sub S$ such that $W=W_J\ct$. Then $W_J$ is finite.
We assume that $w\m w^*$ is the unique automorphism $w\m w^*$ of $W$ such that
$w^*=w_Jww_J$ for $w\in W_J$ and
$t^*=w_Jt\i w_J$ for $t\in\ct$. This automorphism preserves $S$ and has square $1$. We have $J^*=J$
and ${}^JW^J\sub\II_*$ (see \cite{L12, 8.2}).
Let $\ct_{dom}=\{t\in\ct;|w_Jt|=|w_J|+|t|\}=\{t\in\ct;w_Jt\in{}^JW\}$.

(a) If $t\in\ct_{dom}$ we have $|w_Jtw_J|=|t|$; hence $|w_Jt|=|w_Jtw_J|+|w_J|$ and
$w_Jt=(w_Jtw_J)\bul w_J$.
\nl
Indeed, $|w_Jtw_J|=|(t^*)\i|=|t^*|=|t|$.

It is known that

(b) if $t,t'$ are in $\ct_{dom}$ then $tt'\in\ct_{\dom}$ and $|tt'|=|t|+|t'|$ hence
$t\bul t'=tt'$ and $w_J\bul (tt')=w_Jtt'$.
\nl
From (a) we see that $\{w_Jt;t\in \ct_{dom}\}\sub {}^JW^J$; in fact this inclusion is an
equality. For $t\in\ct$ we define $[t]\in{}^JW$ by $W_Jt=W_J[t]$; now $t\m[t]$ is
a bijection
$\ct\lra{}^JW$. Under this bijection and the bijection $\ct_{\dom}\lra{}^JW^J$, $t\lra w_Jt$,
the map ${}^J\p:{}^JW@>>>{}^JW^J$ becomes a map

(c) $\p':\ct@>>>\ct_{dom}$.
\nl
The following result describes explicitly the restriction of $\p'$ to $\ct_{dom}$.

(d) For $t\in\ct_{dom}$ we have ${}^J\p(w_Jt)=w_Jt^2$. Hence $\p'(t)=t^2$.
\nl
Using (a),(b) and the definitions we have
$$\align&{}^J\p(w_Jt)=(w_Jt)\bul(w_Jt)^*{}\i=(w_Jtw_J)\bul w_J\bul w_J\bul(w_Jtw_J)^*{}\i
\\&=(w_Jtw_J)\bul w_J\bul(w_Jt^*{}\i w_J)=(w_Jt)\bul t=w_J\bul t\bul t=w_Jt^2.\endalign$$
This proves (d).

\subhead 2.4\endsubhead
In the setup of 2.3 we assume that $W$ is of affine type
$A_1$. We can assume that $S=\{s_1,s_2\}$ and
$J=\{s_1\}$; now $*$ is the identity map.
We shall write $i_1i_2i_3\do$ instead of $s_{i_1}s_{i_2}s_{i_3}\do$.
Then $A=21\in\ct_{dom}$ and in
fact the elements of $\ct_{dom}$ are precisely the powers
$A^m,m\in\NN$.
The elements of ${}^JW$ are $1t,1t2$ with $t\in\ct_{dom}$.
The elements of ${}^JW^{J^*}$ are $1t$ with $t\in\ct_{dom}$.
If $t=A^m,m\in\NN$, we have

${}^J\p(1t)=1A^{2m}$,

${}^J\p(1t2)=1A^{2m+1}$.

\subhead 2.5\endsubhead
In the setup of 2.3 we assume that $W$ is of affine type
$A_2$. We can assume that $S=\{s_1,s_2,s_3\}$ and
$J=\{s_1,s_2\}$ where $*$
interchanges $s_1,s_2$ and keeps $s_3$ fixed.
We shall write $i_1i_2i_3\do$ instead of $s_{i_1}s_{i_2}s_{i_3}\do$.
Then $A=3121,B=321321,C=312312$ belong to $\ct_{dom}$ and in
fact the elements of $\ct_{dom}$ are precisely the products

(a) $t=A^mB^nC^p$ with $m,n,p$ in $\NN$.
\nl
(Note that $A,B,C$ commute and $A^3=BC$.)
The elements of ${}^JW$ are:

$121t,121t3,121t31,121t32, 121t312,121t321$
\nl
with $t\in\ct_{dom}$.
The elements of ${}^JW^{J^*}$ are $121t$ with $t\in\ct_{dom}$.
If $t$ is as in (a) we have

${}^J\p(121t)=121A^{2m}B^{2n}C^{2p}$,

${}^J\p(121t3)=121A^{2m+1}B^{2n}C^{2p}$,

${}^J\p(121t31)=121A^{2m}B^{2n}C^{2p+1}$,

${}^J\p(121t32)=121A^{2m}B^{2n+1}C^{2p}$,

${}^J\p(121t321)=121A^{2m+2}B^{2n}C^{2p}$,

${}^J\p(121t312)=121A^{2m+2}B^{2n}C^{2p}$.

\enddocument